\documentclass[a4paper,10pt]{article}
\usepackage{amstext}
\usepackage{amssymb,amsmath}
\usepackage{amsfonts}           

\def\R{\mathbb{R}}
\def\CC{\mathbb{C}}

\def\0{\mathbf{0}}
\def\1{\mathbf{1}}
\def\cc{\mathbf{c}}
\def\d{\mathbf{d}}

\def\hu{\hat{\mathbf{u}}}
\def\hv{\hat{\mathbf{v}}}

\def\r{\mathbf{r}}
\def\x{\mathbf{x}}

\def\X{\mathbf{X}}
\def\u{\mathbf{u}}
\def\v{\mathbf{v}}

\def\y{\mathbf{y}}

\def\z{\mathbf{z}}

\def\I{\mathbf{I}}
\def\M{\mathbf{M}}

\def\W{\mathbf{W}}
\def\A{\mathbf{A}}

\def\C{\mathbf{C}}
\def\DD{\mathbf{D}}
\def\VOL{\mathbf{VOL}}
\def\EE{\mathbf{E}}
\def\F{\mathbf{F}}
\def\PP{\mathbf{P}}

\def\RR{\mathbf{R}}
\def\ep{\varepsilon}

\def\diag{\mathrm{diag}}

\def\disc{\mathrm{disc}}

\def\rk{\mathrm{rank}}

\def\Vol{{\mathrm{Vol}}}
\def\VOL{{\mathrm{VOL}}}

\def\part{\cal P} 

\newtheorem{theorem}{Theorem}
\newtheorem{proposition}{Proposition}

\newtheorem{definition}{Definition}
\newtheorem{lemma}{Lemma}

\begin{document}

\begin{center}
\textbf{Relating multiway discrepancy and singular values of graphs and
contingency tables} 
\end{center}

\begin{center}
Marianna~Bolla
\end{center}


\begin{center}
\textit{Institute of Mathematics, Budapest University of Technology and 
Economics} \\ E-mail: {\tt {marib@math.bme.hu}}
\end{center}


\renewcommand\abstractname{Abstract}
                                             
\begin{abstract}

\noindent 
The $k$-way discrepancy $\disc_k (\C)$ of a
rectangular array $\C$ of nonnegative entries is the minimum of 
the maxima of the within- and between-cluster discrepancies
that can be obtained by simultaneous $k$-clusterings (proper partitions) 
of its rows and columns.
In Theorem~\ref{fotetel}, irrespective of the size of $\C$,
we give the following estimate for the
$k$th largest
non-trivial singular value of the normalized table: 
$s_k \le 9\disc_{k } (\C )  (k+2 -9k\ln \disc_{k } (\C ))$, provided
$\disc_{k } (\C ) <1$  and $k\le \rk (\C )$.  
This statement is the converse of Theorem 7 of Bolla~\cite{Bolla14}, 
and the proof
uses some lemmas and ideas of Butler~\cite{Butler}, where only the
$k=1$ case is treated, in which case our upper bound is the tighter. 
The result naturally extends to the singular
values of the normalized adjacency matrix  of a weighted undirected or
directed graph. 



\end{abstract}

\section {Introduction}\label{intro}

In many applications, for example when microarrays are analyzed, our
data are collected in the form of an $m\times n$ rectangular array 
$\C=(c_{ij})$ of
nonnegative real entries, called contingency table. 
We assume that $\C$ is
non-decomposable, i.e., $\C \C^T$ (when $m\le n$) or 
$\C^T \C$  (when $m > n$) is irreducible.
Consequently,
the row-sums
$d_{row,i} =\sum_{j=1}^n c_{ij}$ and column-sums $d_{col,j}=\sum_{i=1}^m c_{ij}$
of $\C$ are strictly positive, and the diagonal matrices
$\DD_{row} =\diag (d_{row,1} ,\dots ,d_{row,m})$ and 
$\DD_{col} =\diag (d_{col,1} ,\dots ,d_{col,n})$ are regular.
Without loss of generality, we also assume that
$\sum_{i=1}^n \sum_{j=1}^m c_{ij} =1$, since neither our main object, the
normalized table 
\begin{equation}\label{cnor}
 \C_{nor} =  \DD_{row}^{-1/2} \C  \DD_{row}^{-1/2} ,
\end{equation}
nor the multiway discrepancies to be introduced are affected by the scaling 
of the entries of $\C$.
It is well known (see e.g.,~\cite{Bolla14}) that the singular values  of 
$\C_{nor}$ are in the [0,1]
interval. Enumerated in non-increasing order, they are the real numbers
$$
 1=s_0 >s_1 \ge \dots \ge s_{r-1} > s_{r} = \dots = s_{n-1} =0  ,
$$
where $r= \rk (\C )$. When $\C$ is non-decomposable,
1 is a single singular value, and it is
denoted by $s_0$, since it belongs to the trivial singular vector pair,
which will be disregarded  in some further calculations. 

Our purpose is to find relations between the $k$th nontrivial singular value
$s_k$ of $\C_{nor}$ and the minimum $k$-way discrepancy of $\C$ defined
herein.

\begin{definition}\label{diszkrepancia}
The multiway discrepancy  of the rectangular array $\C$ of nonnegative entries
in the proper $k$-partition $R_1 ,\dots ,R_k$ of its rows and
$C_1 ,\dots ,C_k$ of its columns is
\begin{equation}\label{disk} 
 \disc (\C ; R_1 ,\dots ,R_k , C_1 ,\dots ,C_k ) =
 \max_{\substack{1\le a\le b\le k \\X\subset R_a , \, Y\subset C_b}} 
 \frac{|c (X, Y)-\rho (R_a,C_b ) \Vol (X)\Vol (Y)|}{\sqrt{\Vol(X)\Vol(Y)}} ,
\end{equation}
where
$c (X, Y) =\sum_{i\in X} \sum _{j\in Y} c_{ij}$ is the cut between
$X\subset R_a$ and $Y\subset C_b$, 
$\Vol (X) = \sum_{i\in X} d_{row,i}$ is the volume of the row-subset $X$, 
$\Vol (Y) = \sum_{j\in Y} d_{col,j}$ is the volume of the column-subset $Y$, 
whereas
$\rho (R_a,C_b) =\frac{c(R_a,C_b)}{ \Vol (R_a) \Vol (C_b)}$ denotes the relative
density between $R_a$ and $C_b$.
The minimum $k$-way discrepancy  of  $\C$ itself is
$$
 \disc_k (\C ) = \min_{\substack{R_1 ,\dots ,R_k \\ C_1 ,\dots ,C_k } } 
 \disc (\C ; R_1 ,\dots ,R_k , C_1 ,\dots ,C_k ).
$$
\end{definition}

In Section~\ref{conc}, I will extend this notion to an edge-weighted graph $G$
and denote it by $\disc_k (G)$. In that setup, $\C$ plays the role of 
the edge-weight matrix
(symmetric in the undirected; quadratic, but usually not symmetric 
in the directed case; and it is the adjacency matrix if $G$ is a simple graph
when the eigenvalues of the normalized adjacency matrix enter into the
estimates, in their decreasing absolute values).
 
Note that $\disc (\C ; R_1 ,\dots ,R_k , C_1 ,\dots ,C_k )$ is the smallest
$\alpha$ such that for every $R_a ,C_b$ pair and for every 
$X\subset R_a$, $Y\subset C_b$,
\begin{equation}\label{dif}
 |c (X, Y)-\rho (R_a,C_b ) \Vol (X)\Vol (Y)| \le \alpha \sqrt{\Vol(X)\Vol(Y)}
\end{equation}
holds.
Hence, in the $k$-partitions of the rows and columns, 
giving the minimum $k$-way discrepancy (say, $\alpha^*$) of $\C$,
every $R_a ,C_b$ pair is $\alpha^*$-regular in terms of the volumes, and
$\alpha^*$ is the smallest possible discrepancy that can be attained
with proper $k$-partitions. 
It resembles the notion of $\epsilon$-regular pairs in the Szemer\'edi
regularity  lemma~\cite{Szemeredi}, albeit with given number of 
vertex-clusters, which are usually not equitable;
further, with volumes, instead of cardinalities. 

Historically,
the notion of discrepancy together with the expander mixing lemma 
was introduced for simple, regular graphs, see 
e.g., Alon, Spencer, Hoory, Linial, Widgerson~\cite{AlonS,Hoory}, 
and extended to Hermitian matrices 
in Bollob\'as, Nikiforov~\cite{BollobasN}.
In Chung, Graham, Wilson~\cite{Chung1}, the authors use the term quasirandom 
for simple
graphs that satisfy any of some equivalent properties, some of them
closely related to discrepancy and eigenvalue separation.
Chung and Graham~\cite{Chung2} prove that for simple graphs `small' discrepancy 
$\disc (G)$ (with our notation, $\disc_1 (G)$) 
is caused by 
eigenvalue `separation': the second largest singular value (which is also
the second largest absolute value eigenvalue), $s_1$,  of the
normalized adjacency matrix is `small', i.e., separated from the
trivial singular value $s_0 =1$,
which is the edge of the spectrum.
More exactly, they prove $\disc (G) \le s_1$, hence giving some kind of
generalization of the expander mixing lemma for \textit{irregular} graphs.

In the other direction, for Hermitian matrices,
Bollob\'as and Nikiforov~\cite{BollobasN} estimate the second largest singular
value of an $n\times n$ Hermitian matrix $\A$ by $C \disc (\A ) \log n$, 
and show that this is best possible up to a multiplicative constant.
Bilu and Linial~\cite{Bilu} prove the  converse of the expander mixing
lemma for simple regular graphs, but their key
Lemma 3.3, producing this statement, goes beyond regular graphs.
In Alon et al.~\cite{Alon10}, 
the authors relax the notion of eigenvalue separation to essential
eigenvalue separation (by introducing a parameter for it, and requiring the
separation only for the eigenvalues of a relatively large part of the graph).
Then they prove relations between the constants of this kind of 
eigenvalue separation and discrepancy. 

For a general rectangular array $\C$ of nonnegative entries, 
Butler~\cite{Butler} proves the following forward 
and backward statement in the $k=1$ case: 
\begin{equation}\label{but}
 \disc (\C ) \le s_1 \le 150\disc (\C ) (1-8\ln \disc (\C ) ) ,
\end{equation}
where $\disc (\C )$ is our $\disc_1 (\C )$ and, with our notation, $s_1$ is 
the largest nontrivial singular value of $\C_{D}$ (he denotes is with 
$\sigma_2$).
Since $s_1 <1$, the upper estimate makes sense for very small discrepancy,
in particular, for
$\disc (\C ) \le 8.868 \times 10^{-5}$. 
The lower estimate further generalizes the expander mixing lemma to
rectangular matrices, but it can be proved with the same tools as in
the quadratic case (see Proposition~\ref{EML} in Section~\ref{conc}).

So far, the overall discrepancy  has been considered in the sense,
that $\disc (\C )$ or $\disc (G)$ measures the largest possible deviation
between the actual and expected connectedness of arbitrary (sometimes disjoint)
subsets $X,Y$, where under expected the hypothesis of 
independence is understood (which corresponds to the rank 1 approximation).  
Note than in~\cite{Butler,Butler1}, $\disc_t (G)$ (or $AltDisc_t (G)$ for
alternating walks in directed graphs) is also introduced, 
which measures the minimum possible deviation between
the actual and expected number of walks of length $t$ between the 
vertex-subsets. Similar notion appears in~\cite{Chung2}, and other notions 
of discrepancy are also introduced in~\cite{Chung3}; for example, the
skew-discrepancy for directed graphs.
Notwithstanding, these papers  consider variants of  the overall 
discrepancy, which corresponds to the one-cluster situation.

My purpose is, in the multicluster scenario, to find 
similar relations between the minimum $k$-way discrepancy and
the SVD of the normalized matrix, for given $k$.
In one direction, in Section~\ref{biz},  I will prove the following. 

\begin{theorem}\label{fotetel}
For every non-decomposable contingency table $\C$
and integer $1\le k\le \rk (\C )$, 
$$
 s_k  \le 9\disc_{k } (\C )  (k+2 -9k\ln \disc_{k } (\C ))  ,
$$
provided $\disc_{k } (\C ) <1$,
where $s_k$  is the $k$th largest non-trivial singular value of the normalized
table $\C_{nor}$ introduced in~(\ref{cnor}).
\end{theorem}

Note that $\disc_k (\C ) =0$ only if $\C$ has a block structure with $k$ 
row- and column-blocks, in which case $s_k =0$ also holds. 
Likewise, $\disc_{k } (\C ) <1$ is not a peculiar requirement, since in view
of $s_k <1$, the upper bound of the theorem has relevance only for 
$\disc_k (\C )$ much smaller than 1; for example, for 
$\disc_{1 } (\C ) \le 1.866\times 10^{-3}$, 
$\disc_{2 } (\C ) \le 8.459\times 10^{-4}$,
$\disc_{3 } (\C ) \le 5.329\times 10^{-4}$, etc.

In the other direction,
in Theorem 7 of~\cite{Bolla14}, I showed that
(under some balancing conditions on the margins and cluster sizes) 
a bit modified version of this
$k$-way discrepancy is $O (\sqrt{2k} S_k +s_k )$, where $S_k$ is the sum of
the squareroots of the $k$-variances of the optimal row- and 
column-representatives (they depend on the normalized singular vectors
corresponding to $s_1 ,\dots ,s_{k-1}$).  In fact, $S_k$ the smaller,
the larger the gap between $s_k$ and $s_{k-1}$ is. 
I will better explain this notion in Section~\ref{last}.
There I will also 
illustrate that $S_k =0$ holds in many special cases, and consequently,
my upper estimate for the $k$-way discrepancy boils down to $B s_k$ with
some absolute constant $B$. 
For example, in the simple graph case,
when $k=2$ and our graph is bipartite, biregular, 
the discrepancy between the two independent vertex-sets is
estimated from above with $B s_2$ by my result, and, up to a constant factor,
this is the same as the estimate proved in Evra et al.~\cite{Evra}. 
In Section~\ref{last}, I will also mention some spectral relations to 
the weak Szemer\'edi regularity lemma~\cite{Borgs,Frieze,Gharan,Szegedy}.

\section{Proof of Theorem~\ref{fotetel}}\label{biz}

Before proving the theorem, I encounter some lemmas of others that I will
use, possibly with some modifications.

Lemma 3 of Bollob\'as and Nikiforov~\cite{BollobasN} is the key to prove their
main result. 
This lemma states that to every $0<\ep <1$ and vector $\x\in \C^n$, 
$\| \x \| =1$, there exists a vector $\y\in \C^n$ such that its coordinates
take no more than 
$\left\lceil \frac{8\pi}{\ep} \right\rceil \left\lceil \frac4{\ep} \log 
\frac{2n}{\ep} \right\rceil$ 
distinct values and $\|\x -\y \| \le \ep$.
This is why $\log n$ appears in their estimate for the second largest singular
value of an $n\times n$ Hermitian matrix. 
Since I do not want to appear the
log-sizes in my estimate in the miniature world of $[0,1]$, 
I will rather use the construction of the following lemma, 
which is indeed a consequence of Lemma 3 of~\cite{BollobasN}.

\begin{lemma}[Lemma 3 of Butler~\cite{Butler}]\label{l1}
To any vector $\x \in \C^n$, $\| \x \| =1$ and diagonal matrix 
$\DD$ of positive real diagonal entries, one can construct a 
step-vector $\y \in \CC^n$ such that $\|\x-\DD\y\|\le \frac13$,
$\| \DD \y \| \le 1$, and the nonzero entries of $\y$ are of the form
$\left( \frac45 \right)^j e^{\frac{\ell }{29} 2\pi i}$ with appropriate
integers $j$ and  $\ell$  ($0\le \ell \le 28$).
\end{lemma}
Note that starting with an $\x$ of real coordinates, we do not need all
the 29 values of $\ell$, only two of them will show up, as it follows from a
better understanding of the construction of~\cite{Butler}. In fact, by
the idea of~\cite{BollobasN}, $j$'s come from dividing the coordinates of
$\DD^{-1} \x / \| \DD^{-1} \x \|$ in decreasing absolute values into groups, 
where
the cut-points are powers of $\frac45$.  
With the notation $\x =(x_s))_{s=1}^n$,  if $x_s$ is
in the $j$-th group, then the corresponding coordinate of the approximating 
complex vector $\y =(y_s )_{s=1}^n$ is as follows. If $x_s =0$, then $y_s=0$, 
otherwise $y_s =\left(\frac45 \right)^j  
e^{\left( \lfloor \frac{29\theta}{2\pi} \rfloor /29 \right) 2\pi i }$, 
where $\theta$ is the argument of $x_s$, $0\le \theta <2\pi$,
and therefore,  $\ell = \lfloor \frac{29\theta}{2\pi} \rfloor$ is an integer
between 0 and 28. However, when the coordinates of $\x$ are real numbers,
then only the values 0 and 14 of $\ell$  can occur, 
since $\theta$ can take only one
of the values 0 or $\pi$, depending on whether $x_s$ is positive or negative.
We will intensively use this observation in our proof.

\begin{lemma}[Lemma 4  of Butler~\cite{Butler}]\label{l2}
Let $\M$ be a matrix with  largest singular value 
$\sigma$ and corresponding unit-norm singular vector pair $\v, \u$. If  
$\x$ and $\y$ are vectors such that $\| \x \|\le 1$, 
$\| \y \|\le 1$, $\| \v -\x \| \le \frac13$,  $\| \u -\y \| \le \frac13$,
then $\sigma \le \frac92 \langle \x , \M \y \rangle $.
\end{lemma}

Note that, in our case, $\M$ is a real matrix and so, $\v ,\u$ have real
coordinates; still, the approximating (step-vectors) $\x , \y$ may have
complex coordinates, and so, 
$\langle .,. \rangle$ denotes the (possibly complex) inner product.
Note that in the possession of real (column) vectors $\x ,\y$ and matrix $\M$,  
$\langle .,. \rangle$ can be written in terms of matrix-vector multiplications
with transpositions:
$\langle \x , \M \y \rangle =\x^T \M \y$. 

\noindent
\textbf{Proof} (of the main theorem).
Assume that $\alpha :=\disc_k (\C ) <1$ and it is attained with the proper 
$k$-partition  
$R_1 ,\dots ,R_k$ of the rows and $C_1 ,\dots ,C_k$ of the columns of $\C$;
i.e., for every $R_a ,C_b$ pair and  
$X\subset R_a$, $Y\subset C_b$ we have
\begin{equation}\label{reg}
| c (X, Y) -\rho (R_a,C_b) \Vol (X) \Vol (Y)| \le \alpha 
 \sqrt{\Vol (X ) \Vol (Y )} .
\end{equation}

Our purpose is to put Inequality~(\ref{reg}) in matrix form by using
indicator vectors and introducing the $m\times n$ auxiliary matrix 
\begin{equation}\label{F}
  \F =\C - \DD_{row} \RR \DD_{col}  ,
\end{equation}
where $\RR =(\rho (R_a,C_b) )$ is the $m\times n$ block-matrix of $k\times k$
blocks with entries equal to $\rho (R_a,C_b)$ over the block $R_a \times C_b$.
With the indicator vectors  $\1_X$ and $\1_Y$ of $X\subset R_a$ and
$Y\subset C_b$, Inequality~(\ref{reg}) has the following equivalent form:
\begin{equation}\label{ind}
 |\langle \1_X , \F \1_Y \rangle  |
 \le \alpha \sqrt{ \langle \1_X ,\C  \1_n \rangle 
   \langle \1_m , \C \1_Y \rangle }
\end{equation}
where $\1_n$ denotes the all 1's vector of size $n$ and
$\langle .,. \rangle$ denotes the (possibly complex) inner product.
Note that in the possession of real (column) vectors and matrices,  
$\langle .,. \rangle$ can be written in terms of matrix-vector multiplications
with transpositions; for example,
$\langle \1_X , \F \1_Y \rangle =\1_X^T \F \1_Y$. 
At the same time, Equation~(\ref{F}) yields 
$$
  \DD_{row}^{-1/2} \F  \DD_{col}^{-1/2} =  \DD_{row}^{-1/2} \C 
 \DD_{col}^{-1/2} -  \DD_{row}^{1/2} \RR \DD_{col}^{1/2}  =
 \C_{nor} - \DD_{row}^{1/2} \RR \DD_{col}^{1/2}  .
$$
Since the rank of the matrix $\DD_{row}^{1/2} \RR \DD_{col}^{1/2} $
 is at most $k$, 
by Theorem 3 of Thompson\footnote{Actually, Thompson stated the
theorem for square matrices, but in the possession of a rectangular one,
we can supplement it with zero rows or columns to make it quadratic; further,
the nonzero singular values of the so obtained square matrix are the same
as those of the rectangular, supplemented with additional zero singular
values that will not alter the shifted interlacing facts.}~\cite{Thompson},
describing the effect of rank $k$ 
perturbations for the singular values, we obtain the following upper estimate
for $s_k$, that is the $(k+1)$th largest (including the trivial 1)
singular value of $\C_{nor}$:
$$
 s_k \le s_{max} (\DD_{row}^{-1/2} \F  \DD_{col}^{-1/2}) =
 \| \DD_{row}^{-1/2} \F  \DD_{col}^{-1/2} \| ,
$$
where $\| .\|$ denotes the spectral norm.

Let $\v \in \R^m$ be the left and $\u\in \R^n$ be the right unit-norm singular 
vector  corresponding to the
maximal singular value of $\DD_{row}^{-1/2} \F  \DD_{col}^{-1/2}$, i.e.,
$$
 |\langle \v , (\DD_{row}^{-1/2} \F  \DD_{col}^{-1/2} ) \u \rangle | = 
 \| \DD_{row}^{-1/2} \F  \DD_{col}^{-1/2} \|.
$$
In view of Lemma~\ref{l1}, 
there are  stepwise constant vectors $\x\in \CC^m$ and $\y \in\CC^n$ such that
$\| \v -\DD_{row}^{1/2} \x \| \le \frac13$ and 
$\| \u -\DD_{col}^{1/2} \y \| \le \frac13$; further,
$\| \DD_{row}^{1/2} \x \| \le 1$ and $\| \DD_{col}^{1/2} \y \| \le 1$.
Then Lemma~\ref{l2} yields
$$
 \| \DD_{row}^{-1/2} \F  \DD_{col}^{-1/2} \| \le
 \frac92 \left| \langle (\DD_{row}^{1/2} \x ), (\DD_{row}^{-1/2} \F  
 \DD_{col}^{-1/2} ) (\DD_{col}^{1/2} \y ) \rangle \right| =
 \frac92 |\langle \x , \F \y \rangle | .
$$
Now we will use the construction in the proof of the Lemma 3~\cite{Butler}
in the special case when the vectors
$\v=(v_s))_{s=1}^m$ and $\u=(u_s))_{s=1}^n$, to be approximated, 
have real coordinates. Therefore, only the following 
three types of coordinates
of the approximating complex vectors $\x =(x_s))_{s=1}^m$ and 
$\y =(y_s )_{s=1}^n$ will appear.
If $v_s =0$, then $x_s=0$ too; if $v_s >0$, then $x_s =(\frac45 )^j$ with
some integer $j$;  if $v_s <0$, then $x_s =(\frac45 )^j 
e^{\frac{28}{29} \pi i}$ with some integer $j$. Likewise,  
if $u_s =0$, then $y_s=0$ too; if $u_s >0$, then $y_s =(\frac45 )^{\ell}$ with
some integer $\ell$;  if $u_s <0$, then $y_s =(\frac45 )^{\ell} 
e^{\frac{28}{29} \pi i}$ with some integer $\ell$.
With these observations, the step-vectors $\x$ and $\y$ can be written as
the following finite sums with respect to the integers $j$ and $\ell$:
$$
 \x = \sum_{j} (\frac45 )^j \x^{(j)} , \quad
 \x^{(j)} =\sum_{a=1}^k ( \1_{{\cal X}_{ja1}} + e^{\frac{28}{29} \pi i} 
  \1_{{\cal X}_{ja2}} ) , \quad \textrm{where}  
$$
$$
 {\cal X}_{ja1} = \{ s: \, x_s  =(\frac45 )^j , \, s\in R_a \}
\quad \textrm{and} \quad
 {\cal X}_{ja2} = \{ s: \, x_s  =(\frac45 )^j e^{\frac{28}{29} \pi i} , 
  \, s\in R_a \} ;
$$ 
likewise,
$$
 \y = \sum_{\ell} (\frac45 )^{\ell} \y^{(\ell )} , \quad
 \y^{(\ell )} =\sum_{b=1}^k ( \1_{{\cal Y}_{\ell b1}} + e^{\frac{28}{29} \pi i} 
  \1_{{\cal Y}_{\ell b2}} ) , \quad \textrm{where}  
$$
$$
 {\cal Y}_{\ell b1} =\{ s: \, y_s =(\frac45 )^{\ell} , \, s\in C_b \}
\quad \textrm{and} \quad  
 {\cal Y}_{\ell b2} = \{ s: \, y_s =(\frac45 )^{\ell} e^{\frac{28}{29} \pi i},
 \, s\in C_b \} . 
$$
Then
\begin{equation}\label{hiv}
\begin{aligned}
 |\langle \x^{(j)}, \F \y^{(\ell )} \rangle | 
 &\le \sum_{a=1}^k  \sum_{b=1}^k \sum_{p=1}^2 \sum_{q=1}^2 
 \left| \langle \1_{{\cal X}_{jap}} ,\F  \1_{{\cal Y}_{\ell bq}} \rangle 
 \right|  \\
 &\overset{(\ref{ind})}\le \sum_{a=1}^k  \sum_{b=1}^k  \sum_{p=1}^2 \sum_{q=1}^2 
\alpha \sqrt{ \langle \1_{{\cal X}_{jap}} ,\C  \1_n \rangle 
              \langle \1_m , \C \1_{{\cal Y}_{\ell bq}} \rangle } \\
 &\le \alpha 2k \sqrt{\sum_{a=1}^k  \sum_{b=1}^k  \sum_{p=1}^2 \sum_{q=1}^2
 \langle \1_{{\cal X}_{jap}} ,\C  \1_n \rangle 
              \langle \1_m , \C \1_{{\cal Y}_{\ell bq}} \rangle }  \\
 &= 2k \alpha \sqrt{ \langle \sum_{a=1}^k \sum_{p=1}^2 
    \1_{{\cal X}_{jap}} ,\C  \1_n \rangle 
 \langle \1_m , \C \sum_{b=1}^k \sum_{q=1}^2 \1_{{\cal Y}_{\ell bq}} \rangle } \\
&= 2k \alpha \sqrt{ \langle |\x^{(j)} | ,\C  \1_n \rangle 
   \langle \1_m , \C |\y^{(\ell )} | \rangle } ,
\end{aligned}
\end{equation}
where in the first inequality we used that $|e^{\frac{28}{29} \pi i}|=1$,
in the second one we used (\ref{ind}), while in the last one,
the Cauchy--Schwarz inequality with $4k^2$ terms. We also introduced the
notation $|\z | = (|z_s |)_{s=1}^n$ for the real vector, the coordinates of
which are the absolute values of the corresponding coordinates of the 
(possibly complex) 
vector $\z$. In the same spirit, let $| \M |$ denote the matrix whose entries
are the absolute values of the corresponding entries of $\M$ (we will use
this only for real matrices). With this formalism, this  is the right moment
to prove the following inequalities that will be used soon to finish the proof:
\begin{equation}\label{kell}
 \sum_{\ell} |\langle \x^{(j)}, \F \y^{(\ell )} \rangle | \le 
  2 \langle | \x^{(j)} | , \C \1_n \rangle , \quad
 \sum_{j} |\langle \x^{(j)}, \F \y^{(\ell )} \rangle | \le 
  2 \langle \1_m , \C  |\y^{(\ell )} | \rangle . 
\end{equation}
Since the two inequalities are of the same flavor, it suffices to prove 
only the first one. Note that it is here, where we use the exact definition of
$\F$ as follows.
$$
\begin{aligned}
\sum_{\ell} |\langle \x^{(j)}, \F \y^{(\ell )} \rangle | 
&\le \langle |\x^{(j)} |, |\F | \sum_{\ell} |\y^{(\ell )} | \rangle   \\
&\le \langle |\x^{(j)} |, (\C + \DD_{row} \RR \DD_{col} ) 
  \1_n \rangle | 
 = 2 \langle | \x^{(j)} | , \C \1_n \rangle  
\end{aligned}
$$
because  $|\y^{(\ell )} |$ is a 0-1 vector and
$\C + \DD_{row} \RR \DD_{col}$ is a (real) matrix of nonnegative entries.
 We also used that the $i$th coordinate of the vector 
$(\C + \DD_{row} \RR \DD_{col} )  \1_n$ for $i\in R_a$ is
$$
 d_{row,i} \left( 1+ \sum_{b=1}^k \rho (R_a ,C_b ) \Vol (C_b ) \right)= 
 2 d_{row,i} 
$$
(here we utilized that the sum of the entries of $\C$ is 1), and therefore,
$$
 (\C + \DD_{row} \RR \DD_{col} )  \1_n = 2 \C \1_n .
$$

Finally, we will finish the proof with similar 
calculations as in~\cite{Butler}. Let us further estimate
$$
 \langle \x , \F \y \rangle = \sum_j \sum_{\ell} 
  \langle (\frac45 )^{j} \x^{(j)} , \F (\frac45 )^{\ell} \y^{(\ell )} 
 \rangle .
$$
Put $\gamma := \log_{4/5} \alpha$; in view of
$\alpha <1$, $\gamma >0$ holds. Then we divide the above summation into
three parts as follows.
$$
\begin{aligned}
 &|\langle \x , \F \y \rangle | \le \sum_j \sum_{\ell} (\frac45 )^{j+\ell}
  |\langle \x^{(j)} , \F  \y^{(\ell )} \rangle | \\ 
&= \underset{\textrm{(a)}}{\sum_{|j-\ell | \le \gamma }} (\frac45 )^{j+\ell}
  |\langle \x^{(j)} , \F  \y^{(\ell )} \rangle | +
 \underset{\textrm{(b)}}{\sum_{j-\ell > \gamma }} (\frac45 )^{j+\ell}
  |\langle \x^{(j)} , \F  \y^{(\ell )} \rangle | +
 \underset{\textrm{(c)}}{\sum_{j-\ell < -\gamma }} (\frac45 )^{j+\ell}
  |\langle \x^{(j)} , \F  \y^{(\ell )} \rangle | .
\end{aligned}
$$
The three terms are estimated separately. Term (a) can be 
bounded from above as follows:
$$
\begin{aligned}
 \sum_{|j-\ell | \le \gamma } (\frac45 )^{j+\ell}
  |\langle \x^{(j)} , \F  \y^{(\ell )} \rangle | & \overset{(\ref{hiv})}{\le}
 2k\alpha   \sum_{|j-\ell | \le \gamma }
  \sqrt{ (\frac45 )^{2j} \langle |\x^{(j)} | ,\C  \1_n \rangle 
   (\frac45 )^{2\ell} \langle \1_m , \C |\y^{(\ell )} | \rangle }   \\
  & \overset{(*)}{\le}
 k\alpha   \sum_{|j-\ell | \le \gamma }
  \left[ (\frac45 )^{2j} \langle |\x^{(j)} | ,\C  \1_n \rangle +
   (\frac45 )^{2\ell} \langle \1_m , \C |\y^{(\ell )} | \rangle \right] \\
   & \overset{(**)}{\le}  k\alpha (2\gamma +1 )  
 \left[ \sum_{j }(\frac45 )^{2j} \langle |\x^{(j)} | ,\C  \1_n \rangle +
 \sum_{\ell } (\frac45 )^{2\ell} \langle \1_m , \C |\y^{(\ell )} | \rangle \right] ,\\
  & \overset{(***)}{\le}  2k\alpha (2\gamma +1 )   ,
\end{aligned}
$$
where in the first inequality, the estimate of~(\ref{hiv}) and in (*),
the geometric-arithmetic mean inequality were used; (**) comes from the
fact that in summation (a), for fixed $j$ or $\ell$, any term can show
up at most  $2\gamma +1$ times, and (***) is due to the easy observation that
\begin{equation}\label{no}
 \sum_{j }(\frac45 )^{2j} \langle |\x^{(j)} | ,\C  \1_n \rangle =
 \| \DD_{row}^{1/2} \x \|^2 \le 1 , \quad 
 \sum_{\ell} (\frac45 )^{2\ell} \langle \1_m , \C |\y^{(\ell )} | \rangle = 
 \| \DD_{col}^{1/2} \y \|^2 \le 1 .
\end{equation}
Terms (b) and (c) are of similar appearance (the role of $j$ and $\ell$ is
symmetric in them), therefore, we will estimate
only (b). Here  $j-\ell > \gamma$, yielding $j +\ell > 2\ell +\gamma$.
Therefore,
$$
\begin{aligned}
 \sum_{j-\ell >  \gamma } (\frac45 )^{j+\ell}
  |\langle \x^{(j)} , \F  \y^{(\ell )} \rangle | 
 &\le \sum_{\ell } (\frac45 )^{2\ell +\gamma} 
 \sum_j  |\langle \x^{(j)} , \F  \y^{(\ell )} \rangle | \\ 
 &\overset{(\ref{kell})}{\le}  \sum_{\ell } (\frac45 )^{2\ell +\gamma} 
 2 \langle \1_m , \C  |\y^{(\ell )} | \rangle \\
 &= 2 (\frac45 )^{\gamma}  \sum_{\ell } (\frac45 )^{2\ell } 
  \langle \1_m , \C  |\y^{(\ell )} | \rangle \overset{(\ref{no})}{\le}  
   2 (\frac45 )^{\gamma}.
\end{aligned}
$$
where, in the second and third inequalities, (\ref{kell}) and (\ref{no})
were used. 
Consequently, (c) can also be estimated from above with 
$2 (\frac45 )^{\gamma}$.

Collecting the so obtained estimates together, we get 
$$
\begin{aligned}
s_k &\le \frac92 |\langle \x , \F \y \rangle | \le \frac92
 \left[ 2k\alpha (2\gamma +1 ) +4 (\frac45 )^{\gamma} \right] =
 9\alpha \left[ 2k\frac{\ln \alpha}{\ln \frac45} +k+2\right] \\
&\le 9\alpha [2k (-4.5)\ln \alpha +  k +2 ] 
  =9\alpha (k+2 -9k\ln\alpha  ) ,
\end{aligned}
$$
that was to be proved. For $k=1$, our upper bound is
tighter than that of~(\ref{but}). 

\section{Some weaker results}\label{last}

Now about our first attempts to prove something like 
Theorem~\ref{fotetel}, because they may be informative for the reader.
\begin{itemize}
\item
First we wanted to use
Lemma 3 of Bollob\'as and Nikiforov~\cite{BollobasN}, 
since, in addition, it specifies the number of
distinct coordinates of the approximating step-vector.
This lemma states that to every $0<\ep <1$ and vector $\x\in \C^n$, 
$\| \x \| =1$, there is a vector $\y\in \C^n$ such that its coordinates
take no more than 
\begin{equation}\label{step}
 \left\lceil \frac{8\pi}{\ep} \right\rceil \left\lceil \frac4{\ep} \log 
\frac{2n}{\ep} \right\rceil 
\end{equation}
values and $\|\x -\y \| \le \ep$.
 
Note that this lemma implies Lemma 3 of Butler~\cite{Butler},
which states that to any unit-norm vector $\x \in \C^n$ and diagonal matrix 
$\DD$ of positive diagonal entries, one can construct a 
step-vector $\y \in \CC^n$ such that $\|\x-\DD\y\|\le \ep$ and
$\| \DD \y \| \le 1$. Even the construction of the two lemmas are similar.

In our case, $\x\in \R^n $ and we need $1/3$ precision. 
Given the diagonal matrix $\DD$ of
positive diagonal entries, we will now construct a 
step-vector $\y$ of complex entries such that $\|\x-\DD\y\|\le 1/3$, by merely 
using Lemma 3 of~\cite{BollobasN}.  First set
$f:=\| \DD^{-1} \x \|$ and $d:=\|\DD\| =\max_i d_i$. Then, by \cite{BollobasN},
to the unit-norm vector $ \DD^{-1} \x /f$ and to $0<\ep <1$ 
there is a step-vector $\y\in \CC^n$, with the same number of different
coordinates as in~(\ref{step}), such that
$$
 \left\| \frac{\DD^{-1} \x}{f} -\y  \right\| \le \ep .
$$
The step-vector $\z =f\y \in \CC^n$, with the same number of different
coordinates as in $\y$, will do for us, since  with an appropriate 
$\ep$ we can reach that $\| \x -\DD \z \| \le \frac13$. Indeed,
$$
 \ep \ge \left\| \frac{\DD^{-1} \x}{f} -\frac{\z }{f} \right\| = \frac1{f}
 \| \DD^{-1} (\x -\DD \z )\| \ge  \frac1{f} \min_i \frac1{d_i}
 \| \x -\DD \z \| = \frac1{fd} \| \x -\DD \z \| .
$$
Therefore,
$$
\| \x -\DD \z \| \le fd\ep =\frac13
$$
holds with $\ep =\frac1{3fd}$ that cannot exceed $\frac13$, since
$fd \ge 1$. This can be seen from the following argument:
$$
1=\|\x\| = \| \DD \DD^{-1} \x \| \le \| \DD \| \cdot \| \DD^{-1} \x \| =df .
$$ 
Eventually, by the construction of~\cite{BollobasN}, 
$|y_j| \le \frac{|x_j |}{d_j f}$, 
$j=1,\dots ,n$. Therefore,
$|z_j | =f|y_j| \le \frac{|x_j |}{d_j}$, and
$|d_j z_j | \le |x_j |$, $\forall j$. Consequently, $\| \DD\z\| \le \| \x \|=1$.

The main implication of this fact is that the maximal number of distinct
coordinates of the step-vector in Lemma 3 of~\cite{Butler} is also of
order $\log n$, and we wanted to make use of this fact in the first
attempts of the proof of some backward statement. For this purpose,
we managed to prove the following lemma, inspired by Lemma 4 
of~\cite{BollobasN}, though, in a more general setup.
We will give the proof too, since it may be of interest for its
own right. 

\begin{lemma}\label{mine}
Let $\C$ be an $m\times n$ matrix of nonnegative real entries
and let the rows and columns have 
positive real weights $d_{r,i}$'s  and $d_{c,j}$'s  (independently of the
entries of $\C$), which are collected in the main diagonals of the
$m\times m$ and $n\times n$
diagonal matrices $\DD_r$ and $\DD_c$, respectively.
Let $R_1 ,\dots ,R_k$ and $C_1 ,\dots ,C_{\ell}$
be proper partitions of the rows and columns; further, $\x \in \CC^m$ and
$\y \in \CC^n$ be stepwise constant vectors having equal coordinates over 
the index sets corresponding to the partition members of 
 $R_1 ,\dots ,R_k$ and $C_1 ,\dots ,C_l$, respectively. 
The $k\times \ell$ real matrix $\C'=(c'_{ab})$ is defined by 
$$
 c'_{ab} := \frac{c(R_a ,C_b )}{\sqrt{\VOL (R_a ) \VOL (C_b )}} ,
 \quad a=1,\dots k; \, b =1 ,\dots ,\ell ,
$$
where $c(R_a ,C_b)$ is the usual cut of $\C$ between $R_a$ and $C_b$, 
whereas $\VOL (R_a ) =\sum_{i\in R_a} d_{r,i}$ and
$\VOL (C_b ) =\sum_{j\in C_b} d_{c,j}$. Then
$$
 | \langle \x ,\C \y \rangle | \le \| \C' \| \cdot \| \DD_r^{1/2} \x \| \cdot
 \| \DD_c^{1/2} \y \| ,
$$
where $\| \C' \|$ denotes the spectral norm, that is the largest singular
value of the real matrix $\C'$, and the squared norm of a complex vector is
the sum of the squares of the absolute values of its coordinates.
\end{lemma}
Note that here the row- and column-weights have nothing to do with the
entries of $\C$, and the volumes are usually not the ones defined in
Section~\ref{intro}; this is why they are denoted by $\VOL$ instead of $\Vol$.

\noindent
\textbf{Proof of Lemma~\ref{mine}}
For the distinct coordinates of $\x$ and $\y$   we introduce
$$
 x_i := \frac{x'_a}{\sqrt{\VOL (R_a )}} \quad \textrm{if} \quad i\in R_a
\quad \textrm{and} \quad
 y_j := \frac{y'_b}{\sqrt{\VOL (C_b )}} \quad \textrm{if} \quad j\in C_b
$$
with $x'_a$ and $y'_b$ that are coordinates of $\x' \in \CC^k$ and 
$\y' \in \CC^l$. Obviously,  $\| \DD_r^{1/2} \x \| =\| \x' \|$ and
$\| \DD_c^{1/2} \y \| =\| \y' \|$.  Then, using $\bar {}$ for the complex
conjugation,
$$
\begin{aligned}
 | \langle \x ,\C \y \rangle | &= 
\left| \sum_{i=1}^m \sum_{j=1}^n x_i {\bar y}_j c_{ij} \right| =
\left| \sum_{a=1}^k\sum_{b=1}^l\frac{x'_a}{\sqrt{\VOL (R_a )}}
\frac{{\bar y}'_b}{\sqrt{\VOL (C_b)}} c(R_a ,C_b ) \right| \\ 
 &= \left| \sum_{a=1}^k \sum_{b=1}^l x'_a {\bar y}'_b c'_{ab} \right| =
 | \langle \x' ,\C' \y' \rangle |  
 \le s_{max} (\C' ) \cdot \| \x' \| \cdot \| \y' \| \\
&= \| \C' \| \cdot \| \DD_r^{1/2} \x \| \cdot \| \DD_c^{1/2} \y \|
\end{aligned}
$$
by the well-known extremal property of the largest singular value,
which finishes the proof.

Using this lemma and the starting steps of the proof of Theorem~\ref{fotetel},
with the matrix $\F$ defined in~(\ref{F}) and the constructed
step-vectors $\x\in \CC^m$, $\y \in \CC^n$, we have
$$
 s_k \le \| \DD_{row}^{-1/2} \F  \DD_{col}^{-1/2} \| \le
 \frac92 |\langle \x , \F \y \rangle | .
$$
We also know from \cite{BollobasN} and the preliminary argument 
that $\x$ takes on at most $r_1 =\Theta (\log m )$,
and  $\y$ takes on at most $r_2 =\Theta (\log n )$
distinct values, which define the proper partitions
$P_{1} ,\dots ,P_{r_1}$ of the rows and
$Q_{1}, \dots ,Q_{r_2}$ of the columns. 
Let us consider the subdivision of them with
respect to $R_1 ,\dots ,R_k$ and $C_1 ,\dots ,C_k$. 
In this way, we obtain the proper partition
$P'_{1} ,\dots ,P'_{\ell_1 }$ of the rows and
$Q'_{1}, \dots ,Q'_{\ell_2 }$ of the columns  with at most
 $\ell_1 = k r_1$ and $\ell_2 =k r_2 $ parts. 

Now, we apply
Lemma~\ref{mine} to the matrix $\F$ and to the step-vectors $\x$ and $\y$, 
which  are also
stepwise constant with respect to  the above partitions. The   row-weights
and column-weights are the $d_{row,i}$'s and $d_{col,j}$'s, respectively.
In view of the lemma, the entries of the $\ell_1 \times \ell_2$ matrix $\F'$
are 
$$
 f'_{ab} := \frac{f(P'_{a} ,Q'_{b} )}{\sqrt{\Vol (P'_{a} ) 
 \Vol (Q'_{b} )}} 
$$
and
$$
 | \langle \x ,\F \y \rangle | \le 
 \| \F' \| \cdot \| \DD_{row}^{1/2} \x \|  \| \DD_{col}^{1/2} \y \| 
 \le  \| \F' \| .
$$ 
But by a well-known linear algebra fact,
$$
 \| \F' \| =s_{max} (\F' ) \le \sqrt{\ell_1 \ell_2 } 
 \max_{a \in \{ 1,\dots ,\ell_1 \} }  \max_{b \in \{ 1,\dots ,\ell_2 \} }
  | f'_{ab} | \le \ell \cdot 
 \disc_{\substack{R_1 ,\dots ,R_k \\ C_1, \dots ,C_k } } (\C ) ,
$$  
where $\ell =\sqrt{\ell_1 \ell_2 } $  and we used Formula (\ref{disk}) 
for the discrepancy. Consequently,
$$
 s_k  \le \frac92  \ell \disc_{k } (\C ) 
$$
follows. The drawback is that the upper bound contains 
$\ell =k\sqrt{r_1 r_2}$ which is of order
$\sqrt{\log m \log n}$. Therefore, we prefer the estimate of 
Theorem~\ref{fotetel} that does not contain the sizes of $\C$.

\item
Another dead-end was the attempt with the following matrix 
$\EE$ instead of $\F$ of~(\ref{F}):
\begin{equation}\label{E}
 \EE  =\C - \DD_{row} {\hat \C } \DD_{col} ,
\end{equation}
where ${\hat \C }= \sum_{i=0}^{k-1} s_i \hv_i \hu_i^T $ 
is an $m\times n$ block-matrix of $k\times k$
blocks with entries equal to ${\hat c}_{ab}$ over the block $R_a \times C_b$.
The vectors  $\hv_i \in \R^m$ and $\hu_i \in \R^n$ are stepwise constant 
over the partitions of $R_1 ,\dots ,R_k$ of the rows and
$C_1 ,\dots ,C_k$ of the columns of $\C$, obtained by spectral clustering 
tools. The vectors  $\hv_i$ and $\hu_i$ themselves were constructed 
via several SVDs in the proof of the forward
statement of~\cite{Bolla14}  so that
 $\DD_{row}^{1/2} \hv_i$ and $\DD_{col}^{1/2} \hu_i$ be `close' to $\v_i$ and
$\u_i$, respectively, for $i=1,\dots ,k-1$ (for $i=0$, they coincide),
where
$\v_i \in \R^m , \u_i \in \R^n$ is the unit-norm singular vector pair 
corresponding to $s_i$ $(i=1,\dots ,r)$. In particular, 
$\v_0 =(\sqrt{d_{row,1}} ,\dots ,
\sqrt{d_{row,m}})^T$ and $\u_0 =(\sqrt{d_{col,1}} ,\dots ,\sqrt{d_{col,n}})^T$.

The point is that the so-called error matrix $\EE$ is close to the matrix
$\DD_{row}^{1/2} (\C_{nor} -\sum_{i=0}^{k-1} s_i \v_i \u_i^T ) \DD_{col}^{1/2}$,
and $ \| \C_{nor} -\sum_{i=0}^{k-1} s_i \v_i \u_i^T \| =s_k$. 
If now $\x\in \CC^m$ and $\y \in \CC^n$ are step-vectors such that 
$\| \DD_{row}^{1/2} \x \| \le 1$, $\| \v_k -\DD_{row}^{1/2} \x \| \le \frac13$ and
$\| \DD_{col}^{1/2} \y \| \le 1$, $\| \u_k -\DD_{col}^{1/2} \y \| \le \frac13$,
then,
$$
 s_{k} \le \frac92  \langle (\DD_{row}^{1/2} \x ),  
( \DD_{row}^{-1/2} \C \DD_{col}^{-1/2} -\sum_{i=0}^{k-1} s_i \v_i \u_i^T ) 
 (\DD_{col}^{1/2} \y ) \rangle . 
$$
Here the upper bound is very close to $\frac92 |\langle \x , \EE \y \rangle |$.
The problem is that  $\langle \1_X , \EE \1_Y \rangle$ cannot be
directly related to the discrepancy, like  $\langle \1_X , \F \1_Y \rangle$. 
However, $\F$ and $\EE$ are very `close' to each other, since comparing
Formulas~(\ref{F}) and~(\ref{E}), the difference between the corresponding
entries of the block-matrices $\RR$ and $\hat \C$ is
$$
 |\rho (R_a,C_b) - {\hat c}_{ab}| =\frac{1}{\Vol (R_a)  \Vol (C_b )}
 \left| \sum_{i\in R_a } \sum_{j \in C_b } \eta_{ij} \right| ,
$$
which is the density of the error matrix $\EE =(\eta_{ij})$ between 
$R_a$ and $C_b$.
If this is small enough, we may expect a finer upper estimate for $s_k$, 
based on $\EE$.
\end{itemize}

\section{Conclusions and applications}\label{conc}

\subsection{Undirected graphs}\label{undir}

The notion of multiway discrepancy naturally extends to edge-weighted graphs. 
A weighted
undirected graph $G=(V,\W )$ is uniquely characterized by its weighted
adjacency matrix $\W$, which is symmetric of nonnegative entries and
zero diagonal.
$\DD = \diag (d_1 ,\dots ,d_n )$ is the diagonal \textit{degree-matrix} 
($d_i =\sum_{j=1}^n w_{ij}$), $\Vol (U) =\sum_{i\in U} d_i$ is the
volume of $U\subset V$, and for simplicity
we assume that $\sum_{i=1}^n d_i =1$; it does not hurt the generality, because
neither the normalized  matrix  $\W_D =\DD^{-1/2} \W \DD^{-1/2}$, 
nor the
multiway discrepancies to be introduced are affected by  the scaling of $\W$.
In case of a simple graph, $\W_D$ is the \textit{normalized adjacency matrix}.
Definition~\ref{diszkrepancia} extends to this case as follows. 

\begin{definition}
The multiway discrepancy  of the undirected,
weighted graph $G=(V,\W )$ in the proper $k$-partition
$V_1 ,\dots ,V_k$  of its vertices is
$$
 \disc (G; V_1 ,\dots ,V_k ) =
 \max_{\substack{1\le a\le b\le k  \\ X\subset V_a , \, Y\subset V_b}} 
 \frac{|w (X, Y)-\rho (V_a,V_b ) \Vol (X)\Vol (Y)|}{\sqrt{\Vol(X)\Vol(Y)}}. 
$$
The minimum $k$-way discrepancy  of the undirected weighted graph $G=(V,\W )$ is
$$
 \disc_k (G) = \min_{V_1 ,\dots ,V_k } \disc (G; V_1 ,\dots ,V_k ) .
$$
\end{definition}
A result, analogous to that of Theorem~\ref{fotetel} can now be formulated in 
terms of the
normalized modularity matrix of $G$, defined in~\cite{Bolla11} as follows.
Denoting by $\d =(d_1 ,\dots ,d_n)^T$ the \textit{degree-vector} (of entries
summing to 1), the so-called \textit{modularity matrix} is 
$\M =\W -\d \d^T$, the $(i,j)$ entry of which just measures the deviation of 
$w_{ij}$ (actual connection of vertices $i$ and $j$) from 
$d_i d_j$ (their connection under independent attachment 
with the vertex-degrees as probabilities). 
With the notation $\sqrt{\d }=(\sqrt{d_1} ,\dots ,\sqrt{d_n})^T$,
the \textit{normalized modularity matrix}  is
$$
 \M_D =  \DD^{-1/2} \M \DD^{-1/2} = \W_D -\sqrt{\d} \sqrt{\d}^T .
$$
The spectrum of $\M_D$ is in the [-1,1] interval, and 0 is always an 
eigenvalue with unit-norm eigenvector $\sqrt {\d}$. All the other eigenvalues
are the same as those of $\W_D$, except the trivial one. Indeed, 
1 is a single eigenvalue of 
$\W_D$ with corresponding unit-norm eigenvector 
$\sqrt {\d}$, provided $\W$ is irreducible.
This becomes a zero eigenvalue of $\M_D$ with the same 
eigenvector. In~\cite{Bolla}, I denoted the eigenvalues 
of $\M_D$ in  decreasing absolute values by 
$|\mu_1 | \ge \dots \ge |\mu_{n-1}| \ge \mu_n =0$. Then the absolute values of 
the eigenvalues of $\W_D$ are $1=\mu_0 \ge 
|\mu_1 | \ge \dots \ge |\mu_{n-1}|$, and they are also the singular values:
$s_k =|\mu_k|$, $k=0,\dots ,n-1$.

\begin{proposition}
Let $G= (V, \W)$ be an edge-weighted, undirected graph. Then 
\begin{equation}\label{enyem}
 |\mu_k |  \le 9\disc_{k } (G )  (k+2 -9k\ln \disc_{k } (G )) ,
\end{equation}
where $\mu_k$ is the $k$-th largest absolute value eigenvalue of
the normalized modularity matrix $\M_D$ $(k=1,\dots ,n-1 )$. 
\end{proposition}

Recall that Bilu and  Linial~\cite{Bilu} prove the following converse of the 
expander mixing lemma for $d$-regular simple graphs on $n$ vertices. 
Assume that for
any disjoint vertex-subsets $S,T$: $\vert e (S,T ) -\frac{|S| |T| d}{n} \vert
\le \alpha \sqrt{|S||T|}$. Then all but the largest adjacency eigenvalue of
$G$ are bounded (in absolute value) by 
$O (\alpha (1+\log \frac{d}{\alpha} ))$.
Note that for a $d$-regular graph the adjacency eigenvalues are $d$
times larger than the normalized adjacency ones, and the deviation
between $e(S,T)$ and the one what is expected in a random $d$-regular graph,
is also proportional to our (1-way) discrepancy in terms of the volumes.
Though  they use disjoint subsets $S,T$, their upper estimate for the 
absolute value 
of the second largest (in absolute value) eigenvalue with the (1-way) 
discrepancy $\alpha$ is $C \alpha (1 -A\log \alpha)$ with some 
absolute constants $A,C$.
Hence, the  upper estimate of~(\ref{but}) or that
of~(\ref{enyem}) in the $k=1$ case are reminiscent of this.  

In the other direction, for the $k=1$ case, a straightforward generalization
of the \textit{expander mixing lemma for irregular graphs} is the following.
\begin{proposition}\label{EML}
$$
\disc (G) =\disc_1 (G) \le \| \M_D \| =s_1 = |\mu_1 |,
$$
where $\| \M_D \|$ is the spectral norm of the normalized modularity matrix 
of $G$.
\end{proposition}
Though, with different notation  (sometimes even a stronger version of it) is
proved in~\cite{Bollabeyond,Butler,Chung2},
we give another short proof here.

\noindent
\textbf{Proof.}
Via  separation  theorems for singular values, $s_1 =|\mu_1 |$ is the
maximum of the bilinear form $\v^T \M_D \u$ over the unit sphere.
Let $X,Y\subset V$ be arbitrary, and 
denote by $\1_X , \1_Y\in \R^n$ the indicator vectors of them. Then
$$
\begin{aligned}
 \| \M_D \| &=\max_{ \| \u\| =\| \v\| =1} |\v^T \M_D \u | \ge
 \left| \left( \frac{\DD^{1/2} \1_X}{\| \DD^{1/2} \1_X \|} \right)^T
 \M_D \left( \frac{\DD^{1/2} \1_Y}{\| \DD^{1/2} \1_Y \|} \right) \right| \\
 &=\frac{|\1_X^T \M \1_Y |}{\| \DD^{1/2} \1_X \|\cdot \| \DD^{1/2} \1_Y \| } =
   \frac{| w (X, Y) - \Vol (X) \Vol (Y)|}{\sqrt{\Vol (X)}\sqrt{\Vol (Y)} } .
\end{aligned}
$$
Taking the maxima on the right-hand side over subsets $X,Y\subset V$, 
the desired relation follows.
Note that the estimate is also valid if we take maxima over disjoint
$X,Y$ pairs only.

For an arbitrary $k$ (between 1 and $\rk \W $), in Theorem 3 of~\cite{Bolla}
we proved that under some balancing conditions for the degrees and the
cluster sizes (when $n\to \infty$), and denoting by $V_1 ,\dots ,V_k $ the
clusters obtained by spectral clustering (see the forthcoming explanation), 
the $(V_a, V_b)$ pairs are
$O (\sqrt{2k} S_k + |\mu_k |)$-volume regular $(a\ne b)$ and 
similar statement holds for the subgraphs induced by  $V_a$'s too. 
In fact, inspired by~\cite{Alon10}, 
there we used a bit different notation and concept of $\alpha$-volume
regular pairs, namely,  for every $X\subseteq V_a$, $Y\subseteq V_b$ 
we required
$$
| w (X, Y) -\rho (V_a ,V_b ) \Vol (X) \Vol (Y)| \le \alpha 
\sqrt{\Vol (V_a) \Vol (V_b)} .
$$
In the above formula, 
the right had side contains the squareroots of the volumes of the
clusters, unlike~(\ref{dif}), which contains  the squareroots of the volumes of 
$X$ and $Y$. However, in the spirit of the 
Szemer\'edi regularity lemma~\cite{Szemeredi},
if we require (\ref{dif}) to hold only for $X,Y$'s satisfying 
$\Vol (X) \ge \ep \Vol (V_i )$, $\Vol (Y) \ge \ep \Vol (V_j )$ with some
fixed $\ep$, then the so modified $k$-way discrepancy, $\disc'_k (G)$, is 
$O (\sqrt{2k} S_k + |\mu_k |)$, and so does $\disc_k (G)$. 
Here the partition $V_1,\dots ,V_k $  is defined so that 
it minimizes
the weighted $k$-variance $S_k^2$ of the vertex representatives
$\r_1 ,\dots ,\r_n \in \R^{k-1}$
obtained as row vectors of the $n\times (k-1)$ matrix of column vectors
 $\DD^{-1/2} \u_i$, 
where $\u_i$ is the unit-norm eigenvector 
corresponding to $\mu_i$  $(i=1,\dots ,k-1 )$. The $k$-variance of the
representatives is defined as
\begin{equation}\label{kszoras}
 {S}_k^2 (\X ) =\min_{(V_1 ,\dots ,V_k )}
\sum_{a=1}^k \sum_{j\in V_a } d_j \| \r_j -{ \cc }_a \|^2 , 
\end{equation}
where ${\cc }_a =\frac1{\Vol (V_a ) } \sum_{j\in V_a } d_j \r_j $ is the
weighted center of cluster $V_a$. 
It is the weighted $k$-means algorithm that gives this minimum, and
the point is that the optimum $S_k$ is just the minimum  distance 
between
the eigensubspace corresponding to $\mu_0 ,\dots \mu_{k-1}$ and  the one
of the suitably transformed step-vectors over the $k$-partitions of $V$.       
In~\cite{Bolla} we also discussed 
that, in view of subspace perturbation theorems, the larger the gap
between $|\mu_{k-1}|$ and $|\mu_k |$, the smaller $S_k$ is. 
So the message is, that here the eigenvectors corresponding to
the largest absolute value eigenvalues have to be used, unlike
usual spectral clustering methods which automatically use the bottom
eigenvalues of the Laplacian or normalized Laplacian matrix
(latter one is just $\I -\W_D $). 
The clusters or cluster-pairs of small discrepancy 
behave like expanders or bipartite expanders. In another context,
they resemble the generalized random or quasirandom graphs  
of~Lov\'asz, S\'os, Simonovits~\cite{LovSos,SimonovitsS}. 

In some
special  cases, $S_k =0$, and then,
$\disc_k (G)\le B |\mu_k | = B s_k$ follows from the above results.  
In particular, $S_k=0$ whenever the vectors 
 $\DD^{-1/2} \u_1, \dots ,\DD^{-1/2} \u_{k-1}$ are step-vectors over the
same proper $k$-partition of the vertices. Some examples:

 \begin{itemize}
\item
If $k=1$, then the unit-norm eigenvector corresponding to $\mu_0 =1$ is 
$\u_0 =\sqrt{\d}$, and $\DD^{-1/2} \u_0 =\1$ is the all 1's vector.
Consequently, the variance of its coordinates is $S_1 =0$.
But in this case, by Proposition~\ref{EML}, we already know that 
$\disc (G)$ can be estimated from above merely by $|\mu_1 | =s_1$.
\item
If $k=2$ and $G$ is bipartite, then $\mu_1 =-1$, $s_1 =1$, and
$S_2^2$, i.e., the 2-variance of the coordinates of the transformed eigenvector
corresponding to $\mu_1$ can be small if $|\mu_2 |$ is separated from
$|\mu_1 |=1$ (see also the bipartite expanders of~\cite{Alon0}).
\item
Let $k=2$ and $G$ be bipartite, biregular  on the independent vertex-subsets 
$V_1 ,V_2$. That is, all the edge-weights within $V_1$ or $V_2$ are zeros,
and the 0-1 weights between vertices of $V_1$ and $V_2$ are such that
$d_i =k_1$ if $i\in V_1$ and $d_i =k_2$ if $i\in V_2$ with the understanding
that $|V_1| k_1 = |V_2| k_2$ (both are the total number of edges in $G$). 
It is easy to see that  the unit-norm eigenvector corresponding to the 
eigenvalue  $\mu_1 =-1$ is $\u_1 =\DD^{1/2} \1_{V_1 } -\DD^{1/2} \1_{V_2 }$, and
$\DD^{-1/2} \u_1 = \1_{V_1 } - \1_{V_2 }$. Therefore, the representatives
of vertices of $V_1$ are all 1's, and those of $V_2$ are $-1$'s, so
$S_2 =0$. Consequently, 
$\disc_2 (G)\le B |\mu_2 |$, with some absolute constant $B$. 
Up to a constant, this was another proof of
Lemma 3.2 of Evra et al.~\cite{Evra}. They call their result expander mixing
lemma for bipartite graphs, and use cardinalities instead of volumes,
but in this special case, these cardinalities are proportional to the volumes 
both within $V_1$ and $V_2$.

\item
Let $G_n$ be a generalized random graph over the symmetric $k\times k$ pattern 
matrix  $\PP =(p_{ab})$, i.e.,  there is a proper $k$-partition, 
$V_1 ,\dots ,V_k$, of its vertices  such that 
 $|V_a |=n_a$ $(a=1,\dots ,k)$, $\sum_{a=1}^k n_a =n$, 
and for any $1\le a\le b\le k$, vertices
$i\in V_a$ and $j\in V_b$ are connected independently, 
with the same probability $p_{ab}$. This is the $k$-cluster generalization of
the classical Erd\H os--R\'enyi random graph, see also~\cite{LovSos} for their
generalized quasirandom counterparts.
 In~\cite{Bolla8} we characterized
the adjacency and normalized Laplacian spectra of such graphs, that
extends to their normalized modularity spectra as follows: both $|\mu_k|=s_k$
and $S_k$ tend to zero almost surely when when $n\to\infty$,
under some balancing conditions for the cluster sizes ($\frac{n_a}{n}\ge c$
with some constant $c$, for $a=1,\dots ,k$).
By our results,  it also holds for the $k$-way discrepancy in the clustering
$V_1 ,\dots ,V_k$. However,
this is not surprising, since this almost sure limit for the
$k$-way discrepancy is easily 
obtained with large deviation principles too, see~\cite{Bolla5}.
\end{itemize}

Summarizing, in the $k=1$ case: 
when the second singular value $|\mu_1 |=s_1$ is small
(much smaller than $s_0 =1$),
then the overall discrepancy is small. But for $k>1$, a small $s_k$ is
necessary, but not sufficient for a small $k$-way discrepancy. In addition,
$S_k$ should be small too. With subspace perturbation theorems, it is small
if $s_k$ is much smaller than $s_{k-1}$. Hence, a gap in the normalized 
modularity
spectrum may be an indication for the number of clusters. 
The two directions together may give a hint about the optimal choice of $k$
if a practitioner wants to find a $k$-clustering of the rows and
columns (or just of the vertices of a graph) with 
small pairwise discrepancies. 
If there not exists a fairly `small'  $k$
with this property, then in the worst case
scenario, the Szemer\'edi regularity lemma~\cite{Szemeredi} with an
enormously large number of clusters (which number only depends on
the maximum pairwise discrepancy to be attained, and does not depend on $n$)
comes into existence. Weak versions of this lemma
(where $V_1 , \dots ,V_k $ are not necessarily equitable) are also
available, see e.g.,~\cite{Borgs,LovSzeg}. 


Note that $\M_D$ corresponds to the compact operator taking conditional
expectation between the margins with respect to the symmetric joint 
distribution embodied by $\W$. In~\cite{Bolla} we proved that for given $k$,
the eigenvalues 
$\mu_1 ,\dots ,\mu_{k-1}$ and the corresponding eigensubspace are testable,
consequently $S_k$ is also testable, in the sense of~\cite{Borgs}. This is
important when we have a very large network and want to estimate these 
quantities based on a smaller sample selected with an appropriate randomization
from the large one. We also remark that spectral or operator proofs of the
regularity lemma, together with low-rank constructions, are at our disposal, 
for example,~\cite{Frieze,Gharan,Szegedy}.

\subsection{Directed graphs}\label{dir}

A directed weighted graph $G=(V,\W )$ is described by its quadratic, but 
usually not symmetric
weight matrix $\W=(w_{ij})$ of zero diagonal, 
where $w_{ij}$ is the nonnegative weight of the $i\to j$ edge $(i \ne j )$.
The row-sums
$d_{out,i} =\sum_{j=1}^n w_{ij}$ and column-sums $d_{in,j}=\sum_{i=1}^n w_{ij}$
of $\W$ are the out- and in-degrees, while 
$\DD_{out} =\diag (d_{out,1} ,\dots ,d_{out,n})$ and 
$\DD_{in} =\diag (d_{in,1} ,\dots ,d_{in,n})$ are the diagonal
out- and in-degree matrices, respectively. 
Now Definition~\ref{diszkrepancia} can be formulated as follows. 
\begin{definition}
The multiway discrepancy  of the directed,
weighted graph $G=(V,\W )$ in the in-clustering
$V_{in,1} ,\dots ,V_{in,k}$ and out-clustering $V_{out,1} ,\dots ,V_{out,k}$ 
of its vertices is
$$
\begin{aligned}
 &\disc (G; V_{in,1} ,\dots ,V_{in,k}, V_{out,1} ,\dots ,V_{out,k})  \\
 &=\max_{\substack{1\le a\le b\le k \\ X\subset V_{out,a} , \, Y\subset V_{in,b}}} 
 \frac{|w (X, Y)-\rho (V_{out,a},V_{in,b} ) \Vol_{out} (X)\Vol_{in} (Y)|} 
 {\sqrt{\Vol_{out}(X)\Vol_{in} (Y)}},
\end{aligned}
$$
where $w(X,Y)$ is the sum of the weights of the $X\to Y$ edges, whereas
$\Vol_{out} (X) =\sum_{i\in X} d_{out,i}$ and 
$\Vol_{in} (Y) =\sum_{j\in Y} d_{in,j}$ are the out- and in-volumes,
respectively. 
The minimum $k$-way discrepancy  of the directed weighted graph $G=(V,\W )$ is
$$
 \disc_k (G) = \min_{\substack{V_{in,1} ,\dots ,V_{in,k} \\ 
  V_{out,1} ,\dots ,V_{out,k} }}
 \disc (G; V_{in,1} ,\dots ,V_{in,k}, V_{out,1} ,\dots ,V_{out,k}).
$$
\end{definition}
Butler~\cite{Butler} treats the $k=1$ case, and for a general $k$, 
Theorem~\ref{fotetel} implies the following.
\begin{proposition}
Let $G= (V, \W)$ be directed edge-weighted graph. Then 
$$
 s_k  \le 9\disc_{k } (G )  (k+2 -9k\ln \disc_{k } (G )) ,
$$
where $s_k$ is the $k$-th largest nontrivial singular value  of
the normalized  edge-weight matrix $\W_{D}=\DD_{out}^{-1/2} \W \DD_{in}^{-1/2}$.
\end{proposition}
We applied the SVD based algorithm to find migration patterns in
the set of 75 countries, and found 3 underlying  
immigration and emigrationin trait clusters.
Since the  algorithm is the same as for rectangular
matrices, I will describe it in the next subsection.

\subsection{Back to rectangular arrays}

In multivariate statistics, sometimes our data are collected in an
$m\times n$ matrix $\C$, where the entries are 
frequency  counts 
corresponding to the joint distribution of two categorized random variables
(taking on $m$ and $n$ discrete values, respectively). 
Such a $\C$ is called contingency table in statistical language, and the
data are popularly said to be cross-tabulated.  
The $\chi^2$ statistic, which measures the deviation from independence,
is $N\sum_{i=1}^{r-1}s_i^2$ with my notation, where $N$ is the (usually
`large')  sample size, but
the second factor can be `small' if $s_1$ is `small', and this corresponds to
the existence of a good rank 1 approximation of $\C$. This fact is
also supported by the $\disc (\C) =\disc_1 (\C) \le s_1$ relation.  
Otherwise, one may ask, whether there exists a `good'
rank $k$ approximation for some integer $1<k<r=\rk (\C)$, which problem is
treated in correspondence analysis by the first $k$ dyads of the SVD of
$\C_D$. However, there it is not made exact how $s_k$ is estimated by
$\disc_k (\C )$. Our Theorem~\ref{fotetel} says that if the minimum $k$-way
discrepancy is very `small', i.e., the sub-tables $R_a \times C_b$ behave
like independent tables in the optimal $k$-partitions of the rows and columns,
then $s_k$ is small too.

In the other direction, in~\cite{Bolla14}, we
proved the following. Given the $m\times n$ contingency
table $\C$, consider the spectral clusters $R_1 ,\dots ,R_k$ of its rows and 
$C_1 ,\dots ,C_k$ of its columns,
obtained by applying the $k$-means algorithm for the $(k-1)$-dimensional 
row- and column representatives, defined as the row vectors of the matrices
of column vectors $(\DD_{row}^{-1/2} \v_1 ,\dots ,\DD_{row}^{-1/2} \v_{k-1})$
and $(\DD_{col}^{-1/2} \u_1 ,\dots ,\DD_{col}^{-1/2} \u_{k-1})$, respectively,
where $\v_i , \u_i$ is the unit norm singular vector pair corresponding to $s_i$
$(i=1,\dots ,k-1)$.
In fact, these partitions minimize
the weighted $k$-variances  $S_{k,row}^2$ and 
$S_{k,col}^2$  of these row- and column-representatives (see (\ref{kszoras})).
Then, under some balancing conditions for the margins 
and for the cluster sizes, we proved that 
$\disc_k (\C ) \le B (\sqrt{2k} (S_{k,row} +S_{k,col}) +s_k )$,
with some absolute constant $B$.
This is the base of our algorithm, with fixed $k$.

We remark that the correspondence analysis uses the above $(k-1)$-dimensional 
row- and column-representatives for simultaneously plotting the row- and
column-categories in $\R^{k-1}$ ($k=2,3$ or 4 in most applications), 
and hence, the practitioner can draw
conclusions from their mutual positions.
For example, in microarray analysis we can plot the genes and conditions
together, and the biclusters obtained by $k$-clustering the row- and
column-representatives give clusters of the genes and the conditions
such that, every gene-cluster and condition-cluster pair behaves like
a random weighted bipartite graph in the sense, that genes and conditions of
the same cluster nearly independently influence each other, which fact may have
importance for practitioners. In~\cite{Bolla14} it is also shown that when these
$k$-variances are very `small', then our construction (described there 
with the modified
dyads) for the rank $k$ approximation produces a table of nonnegative entries.
On the contrary, 
a drawback of correspondence analysis is that the automatic low-rank
approximation of the table usually contains negative entries.  

In the possession of networks or microarrays, practitioners
want to find a fairly small $k$, such that there is a $k$-cluster structure
behind the table or the graph in the sense that the subgraphs and bipartite 
subgraphs have  `small' discrepancy. It depends on the table or the graph that 
how small discrepancy can be attained and with what $k$. The above theory
tells that we have to inspect the normalized spectra, together with 
spectral subspaces, since the leading ones carry a lot of information
about the smallest attainable discrepancy. 

\section*{Acknowledgement}

The author wishes to thank Gergely Kiss and Zolt\'an Mikl\'os S\'andor for
discussions on the topic.
Parts of the research were done under the auspices of the 
Budapest Semesters of Mathematics program,
 in the framework of an undergraduate
research course on spectral clustering 
 with the participation of US students James Drain, Cristina Mata, 
Matthew Willian, and in particular, Calvin Cheng whose 
computer processing of real-word data helped in formulating the main theorem.
The research was also supported by the
T\'AMOP-4.2.2.C-11/1/KONV-2012-0001 project.


\begin{thebibliography}{1}


\bibitem{Alon0}
 Alon, N., 1986 Eigenvalues and expanders, \textit{Combinatorica}
\textbf{6} (1986), 
83-96.

\bibitem{AlonS}
Alon, N. and Spencer, J. H., \textit{The Probabilistic Method}, Wiley (2000).

\bibitem{Alon10} 
Alon, N., Coja-Oghlan, A., Han, H., Kang, M., R\"odl, V. and Schacht, M.,
Quasi-randomness and algorithmic regularity for graphs with general
degree distributions, \textit{Siam J. Comput.} \textbf{39} (2010), 2336--2362.

\bibitem{Bilu}
Bilu, Y. and Linial, N., Lifts, discrepancy and nearly optimal spectral gap,
\textit{Combinatorica} \textbf{26} (2006), 495--519. 

\bibitem{Bolla5}
Bolla, M., Recognizing linear structure in noisy matrices,
\textit{Linear Algebra and Its Applications} \textbf{402} (2005), 228-244.

\bibitem{Bolla8}
Bolla, M., Noisy random graphs and their Laplacians,
\textit{Discrete Mathematics} \textbf{308} (2008), 4221-4230.

\bibitem{Bollabeyond}
Bolla, M., Beyond the expanders, \textit{International Journal of 
Combinatorics}, Paper 787596 (2011). 

\bibitem{Bolla11}
Bolla, M., Penalized versions of the Newman--Girvan modularity and their
relation to normalized cuts and k-means clustering,
\textit{Physical Review E} \textbf{84} (1), 016108 (2011). 


\bibitem{Bolla}
Bolla, M., Modularity spectra, eigen-subspaces and structure of weighted
graphs, \textit{European Journal of Combinatorics} \textbf{35} (2014), 
105--116. 

\bibitem{Bolla14} 
Bolla, M., SVD, discrepancy, and regular structure of contingency tables,
\textit{Discrete Applied Mathematics} \textbf{176} (2014), 3-11.

\bibitem{BollobasN}
Bollob\'as, B. and Nikiforov, V.,Hermitian matrices and graphs: singular
values and discrepancy, \textit{Discret. Math.} \textbf{285} (2004), 17--32.

\bibitem{Borgs} 
Borgs, C, Chayes, J. T., Lov\'asz, L., T.-S\'os, V. and 
Vesztergombi, K., Convergent graph sequences I: Subgraph
Frequencies, metric properties, and testing, \textit{Advances in Math.}
\textbf{219} (2008), 1801--1851.

\bibitem{Butler}
Butler, S., Using discrepancy to control singular values for nonnegative
matrices, \textit{Linear Algebra Appl.} \textbf{419} (2006), 486--493.  

\bibitem{Butler1}
Butler, S., Relating singular values and discrepancy of weighted directed 
graphs. In \textit{Proc. of the 17th Annual ACM-SIAM Symposium on Discrete
Algorithms, Miami, FL, 2006}, SIAM, Philadelphia, PA, pp.~1112--1116 (2006).

\bibitem{Chung1} 
Chung, F., Graham, R. and Wilson, R. K., Quasi-random graphs,
\textit{Combinatorica} \textbf{9} (1989), 345--362.

\bibitem{Chung2} 
Chung, F. and Graham, R., Quasi-random graphs with given degree
sequences, \textit{Random Struct. Algorithms} \textbf{12} (2008), 1--19.

\bibitem{Chung3}
Chung, F. and Kenter, F., Discrepancy inequalities for directed graphs,
\textit{Discrete Applied Mathematics} \textbf{176} (2014), 30-42.


\bibitem{Evra}
Evra, S., Golubev, K., Lubotzky, A., Mixing properties and the chromatic
number of Ramanujan complexes, arXiv:1407.7700 [math.CO] (2014).

\bibitem{Frieze} 
Frieze, A. and Kannan, R., Quick approximation to matrices
and applications, \textit{Combinatorica} \textbf{19} (1999), 
175--220. 

\bibitem{Gharan}
Gharan, S. H., Trevisan, L., A new regularity lemma and faster approximation
algorithms for low threshold rank graphs. In. \textit{Proc. APPROX-RANDOM'13},
pp. 303-316 (2013). 

\bibitem{Hoory} 
Hoory, S., Linial, N. and Widgerson, A., Expander graphs
and their applications, \textit{Bull. Amer. Math. Soc. (N. S.)}
\textbf{43} (2006), 439--561.

\bibitem{LovSzeg} 
Lov\'asz, L. and Szegedy, B., 
Szemer\'edi's Lemma for the analyst.
\textit{Geom. Func. Anal.} \textbf{17} (2007), 252--270.

\bibitem{LovSos} 
Lov\'asz, L. and T.-S\'os V., Generalized quasirandom
graphs, \textit{J. Comb. Theory B} \textbf{98} (2008), 146--163.

\bibitem{SimonovitsS} 
Simonovits, M. and  T.-S\'os, V., Szemer\'edi's partition and
quasi-randomness, \textit{Random Struct. Algorithms} \textbf{2} (1991), 1--10. 

\bibitem{Szegedy}
Szegedy, B., Limits of kernel operators and the spectral regularity lemma,
\textit{European J. Combin.} \textbf{32} (2011), 
1156-1167.

\bibitem{Szemeredi}
Szemer\'edi, E., Regular partitions of graphs. In \textit{Colloque
Inter. CNRS. No. 260, Probl\'emes Combinatoires et
Th\'eorie Graphes} (Bermond J-C, Fournier J-C, Las Vergnas M and
Sotteau D eds), pp.~399--401 (1976).

\bibitem{Thompson}
Thompson, R. C., The behavior of eigenvalues and singular values under 
perturbations of restricted rank,
\textit{Linear Algebra Appl.} \textbf{13} (1976), 69--78.  

\end{thebibliography}
\end{document}